\magnification=1200

{\centerline{\bf Small Divisor problems and $A_p$ weights with an application}}

\bigskip
{\centerline {Sagun Chanillo}}
\bigskip
{\centerline{\it To Bernard Helffer in friendship and admiration}}
\bigskip{\centerline{\bf To appear in Pure and Applied Functional Analysis(PAFA)}}

{\centerline{\bf in a special volume for B. Helffer}}

\bigskip\noindent
{\bf Abstract:} We show that the Muckenhoupt theory of weight functions can be used to understand certain small divisor problems. We apply this link to obtain a quantitative version of the Ehrenpreis-Malgrange local solvability theorem for constant coefficient linear differential operators, where the loss of derivatives as measured in the scale of Sobolev spaces is related to the $A_p$ class of the symbol of the differential operator.



\bigskip\noindent
\S 1. {\bf Introduction}

\bigskip
The Muckenhoupt theory of $A_p$ weights has played a fundamental role in Harmonic Analysis since the fundamental paper of Muckenhoupt [10], who established the necessity and sufficiency of weight functions satisfying the $A_p$ condition for weighted norm inequalities to hold for the Hardy-Littlewood maximal function. Since then, $A_p$ weights have proved useful in the study of Potential theory for Harmonic functions on Lipschitz domains, unique continuation problems, Harnack inequalities and L. Caffarelli's work on free boundaries.

Our aim here is to point out a link between small divisor problems and $A_p$ weights. We shall exploit this link by giving a quantitative proof of the Ehrenpreis-Malgrange theorem [5], [6], [7], [9] on local solvability for constant coefficient differential operators and also relate the inequalities we obtain to a certain $A_p$ class to which the symbol of the differential operator belongs to. We have nothing to say if such links are also available for the variable coefficient case investigated by Nirenberg-Treves [12], Beals-C. Fefferman[1] or the over-determined situation studied in Chanillo-Treves[2] and Cordaro-Hounie[4].

We introduce some notation. We set
              $$D=\big( {{1}\over{i}}{{\partial}\over{\partial x_1}}, \cdots, {{1}\over{i}}{{\partial}\over{\partial x_n}}\big),\ i=\sqrt{-1}.$$

We shall prove our theorem in the periodic setting and so we will need the $n$-torus $T^n$ and the lattice points $Z^n$ in $R^n$, where $Z$ denotes the integers. Sobolev spaces $H^s(T^n)$ are defined on this torus via Fourier coefficients $\hat{f}(m)$ by
$$H^s(T^n)=\{f\ |\ \ \sum_{m\in Z^n}|\hat{f}(m)|^2(1+|m|^2)^s<\infty\},$$
with associated norm
            $$||f||_s=\big(\sum_{m\in Z^n}|\hat{f}(m)|^2(1+|m|^2)^s\big)^{{{1}\over{2}}}.$$

We will defer the definition of Muckenhoupt's $A_p$ class to the next section. For the statement of our theorem we will need to introduce $Q_0$, the cube of unit edge length, centered at the origin.

\bigskip\noindent
{\bf Theorem 1:} Let $P(D)$ denote any linear, constant coefficient differential operator. Let $\rho\in R$. There is a subset $S\subseteq Q_0$ of full measure, that is $|S|=|Q_0|$, such that for any $\xi_0\in S$, the conjugated equation
             $$ e^{-i(x,\xi_0)}P(D)(e^{i(x,\xi_0)}u)=f$$
has a solution on the torus $T^n$, with the accompanying estimate,
                           $$||u||_\rho\leq C ||f||_{\rho+s},$$
with $s=np>0$. The $p$ is the $A_p$ class to which the polynomial $|P(\xi)|$ belongs to, where $P(\xi)$ is the symbol of the differential operator $P(D)$. The constant $C$ in the inequality is an absolute constant which is independent of $u,f$.

\bigskip

\bigskip\noindent
{\bf Corollary 1:} Let $P(D)$ denote any linear, constant coefficient differential operator. Let $\rho\in R$. Assume $f$ is compactly supported in $T^n$. Then the equation $P(D)u=f$ has a solution $u$ on the torus $T^n$ such that $u=e^{i(x,\xi_0)}U$.
            $$||U||_\rho\leq C ||f||_{\rho+np},$$
and $p$ is as above and $C$ independent of $U,f$.
\bigskip\noindent

Using the fact that $e^{i(x,\xi_0)}$ is a smooth function, and $|e^{i(x,\xi_0)}|=1$, the corollary that follows is immediate from Corollary 1.

\bigskip\noindent
{\bf Corollary 2:} Let $\Omega\subset \Omega^\prime\subset\subset T^n$ be open sets. Let $\eta$ be a smooth cut-off function such that $\eta\equiv 1$ on $\overline{\Omega}$ and supported in $\Omega^\prime$. Then for $f$ compactly supported, there is a solution to $P(D)u=f$ in $T^n$ such that
                     $$||\eta u||_\rho\leq C(\Omega^\prime, n, \rho)||f||_{\rho+np},\ \rho\in R.$$

\bigskip\noindent
We wish to thank Paulo Cordaro for several helpful remarks pertaining to Corollaries 1 and 2. We also thank the anonymous referee for comments that helped improve the paper.

\bigskip\bigskip\noindent
\S 2.{\bf Proof of the theorem}
\bigskip

We will consider here $w(x)>0$, and $w\in L^1_{loc}(R^n)$ called
a weight or a weight function.
\bigskip
\noindent {\bf Definition:} We say $w\in D_p$, doubling of order $p$,
if and only if there exists an absolute constant $C>0$, such that for
any ball $B_s(x_0)$ of radius $s$ centered at $x_0$, we have for all
balls $B_t\subseteq B_s(x_0)$ the inequality below:
         $$\int_{B_s(x_0)} w(x) dx\leq C \biggr({{s}\over{t}}\biggr)^{np}\int_{B_t}
         w(x)\ dx.$$
The ball $B_t$ need not be centered at $x_0$ but it has radius $t$
and obviously since $B_t\subseteq B_s(x_0)$ we must have $t\leq s$.
The constant $C$ is independent of the center $x_0$ and $t,s$ in the
definition. Furthermore, the doubling condition does not depend if we use balls or cubes in the definition.

\bigskip
\noindent {\bf Definition:} Let $w>0$ and locally integrable be a
weight function. For $1<p<\infty$, we say $w\in A_p$, or an
$A_p$ weight, if and only if for all balls $B\subset R^n$, we have
    $$\sup_B \bigg( {{1}\over{|B|}}\int_B w(x)\,
    dx\bigg)\bigg({{1}\over{|B|}}\int_B w^{-{{1}\over{p-1}}}\,
    dx\bigg)^{p-1}\leq C<\infty.$$
When $p=1$ we say $w\in A_1$ if and only if
        $${{1}\over{|B|}}\int_B w(x)\, dx\leq C {\rm ess\ inf}_B
        w(x).$$

\bigskip
\noindent {\bf Proposition 1:} If $w\in A_p$ then $w\in D_p$.

\bigskip
\noindent {\bf Remark 1:} The converse is false. C. Fefferman and Muckenhoupt [8] have constructed
weight functions that are doubling, but lie in no $A_p$ class.

\medskip
\noindent{\bf Proof:} First we have

   $$1={{1}\over{|B|}}\int_B w^{{{1}\over{p}}}w^{-{{1}\over{p}}}$$
Now by Holder's inequality applied to the integral above we have
with ${{1}\over{p}}+{{1}\over{q}}=1$,
$$1\leq \biggr( {{1}\over{|B|}}\int_B w\biggr)^{1/p}\biggr(
{{1}\over{|B|}}\int_B w^{-{{1}\over{p-1}}}\biggr)^{1/q}.$$ We also
the used the fact that ${{p}\over{q}}=p-1$ and so
${{q}\over{p}}=1/(p-1).$

Raising both sides to exponent $p$ and using ${{p}\over{q}}=p-1$ we
get
 $$1\leq \biggr( {{1}\over{|B|}}\int_B w\biggr)\biggr(
{{1}\over{|B|}}\int_B w^{-{{1}\over{p-1}}}\biggr)^{p-1}.\eqno (1)$$

Notice the quantity on the right side above in $(1)$ is exactly what
appears in the definition of $A_p$. In fact we just checked that
$(1)$ holds for any non-negative function.

Now we are ready to prove the proposition. We will suppress the
dependence of the center of the ball so $B_s(x_0)$ will be simply
denoted by $B_s$. We have
  $$\int_{B_s}w(x)\ dx = {{\int_{B_s}w(x)\ dx}\over{{{1}\over{t^n}}\int_{B_t}w(x)\ dx}}\ {{1}\over{t^n}}\int_{B_t}w(x)\ dx.\eqno (2)$$
  Next we re-write the right side of $(2)$ to get,
     $$\biggr({{s}\over{t}}\biggr)^n {\biggr({{1}\over{s^n}}{\int_{B_s}w(x)\ dx\biggr)\biggr({{1}\over{s^n}}\int_{B_s}w^{-{{1}\over{p-1}}}\biggr)^{p-1}}
     \over{\biggr({{1}\over{t^n}}\int_{B_t}w(x)\ dx\biggr)\biggr({{1}\over{s^n}}\int_{B_s}w^{-{{1}\over{p-1}}}\biggr)^{p-1} }}\ \int_{B_t}w(x)\ dx$$

Since by hypothesis $w\in A_p$, the numerator in the expression above is
bounded by a uniform constant independent of the ball $B_s$ and so
the expression above is bounded by
                $$\biggr({{s}\over{t}}\biggr)^n {{C}\over{\biggr({{1}\over{t^n}}\int_{B_t}w(x)\ dx\biggr)
                \biggr({{1}\over{s^n}}\int_{B_s}w^{-{{1}\over{p-1}}}\biggr)^{p-1} }}\ \int_{B_t}w(x)\ dx\eqno (3)$$
         Now since $B_t\subseteq B_s$ we have

         $$\int_{B_s}w^{-{{1}\over{p-1}}}\geq
         \int_{B_t}w^{-{{1}\over{p-1}}}.$$
Hence it follows from the last inequality that,
          $$\biggr({{1}\over{s^n}}\int_{B_s}w^{-{{1}\over{p-1}}}\biggr)^{p-1}\geq
          \biggr({{t}\over{s}}\biggr)^{n(p-1)}\biggr({{1}\over{t^n}}\int_{B_t}w^{-{{1}\over{p-1}}}\biggr)^{p-1}\eqno (4)$$
          Using $(4)$ in the denominator of $(3)$, we get a lower
          bound for the denominator as follows,
$$\biggr({{1}\over{t^n}}\int_{B_t}w(x)\ dx\biggr)
                \biggr({{1}\over{s^n}}\int_{B_s}w^{-{{1}\over{p-1}}}\biggr)^{p-1}
                \geq$$
      $$\biggr({{t}\over{s}}\biggr)^{n(p-1)}\biggr({{1}\over{t^n}}\int_{B_t}w(x)\ dx\biggr)
                \biggr({{1}\over{t^n}}\int_{B_t}w^{-{{1}\over{p-1}}}\biggr)^{p-1}.$$
But by the observation encoded in $(1)$, the expression on the right
above is further lower bounded by
                                   $$\biggr({{t}\over{s}}\biggr)^{n(p-1)}.\eqno (5)$$
Using the lower bound $(5)$ in the denominator of $(3)$ we finally
arrive at,
         $$\int_{B_s}w(x)\ dx \leq
         C\biggr({{s}\over{t}}\biggr)^{n(p-1)}\biggr({{s}\over{t}}\biggr)^n\int_{B_t}w(x)\
         dx.$$
Thus we have verified the $D_p$ condition.

An important consequence of Proposition 1 proved above is

\bigskip
\noindent {\bf Remark 2:} The definition for $A_p$ does not depend on the fact that we use balls. In view of the Doubling property proved above, we get an equivalent definition if we had used cubes.

The relation between different $A_p$ classes is encoded by the two propositions stated and proved below.

\bigskip
\noindent {\bf Proposition 2:} Let $w\in A_p$ then $w\in A_q$ for
all $q\geq p$.

\medskip
\noindent {\bf Proof:} This is a consequence of H\"older's
inequality. First notice that we have by hypothesis
${{1}\over{q-1}}<{{1}\over{p-1}}$. Set $r={{q-1}\over{p-1}}>1$ and
by H\"older, (or averages are monotonically increasing)
        $$\big({{1}\over{|B|}}\int_B w^{-{{1}\over{q-1}}}(x)\,
        dx\big)^{q-1}\leq \big({{1}\over{|B|}}\int_B w^{-{{r}\over{q-1}}}(x)\,
        dx\big)^{{{q-1}\over {r}}}.$$
This yields,
        $$\big({{1}\over{|B|}}\int_B w^{-{{1}\over{q-1}}}(x)\,
        dx\big)^{q-1}\leq \big({{1}\over{|B|}}\int_B w^{-{{1}\over{p-1}}}(x)\,
        dx\big)^{p-1}.$$
Multiplying both sides of the inequality above by
${{1}\over{|B|}}\int_B w(x)\, dx$ we easily see the conclusion of
our Proposition.
\bigskip
\noindent {\bf Proposition 3:} Let $w\in A_p$, then
$w^{-{{1}\over{p-1}}}$ is in $A_{p^\prime}$, where
${{1}\over{p}}+{{1}\over{p^\prime}}=1$ and conversely.

\medskip
\noindent {\bf Proof:} It is enough to prove one direction. By
definition we need to consider
   $$\big({{1}\over{|B|}}\int_B w^{-{{1}\over{p-1}}}(x)\,
   dx\big)\big({{1}\over{|B|}}\int_Bw^{{{1}\over{(p-1)(p^\prime-1)}}}(x)\,
   dx\big)^{p^\prime-1}.$$
But $(p-1)(p^\prime-1)=1$. Using this fact, the expression above is
$$\bigg(\bigg({{1}\over{|B|}}\int_B w(x)\,
    dx\bigg)\bigg({{1}\over{|B|}}\int_B w^{-{{1}\over{p-1}}}\,
    dx\bigg)^{p-1}\bigg)^{p^\prime-1}\leq C.$$

 We next need two theorems whose proof relies on the Calderon-Zygmund cube bisection lemma.

\bigskip \noindent {\bf Theorem 2([3], [10]):} Let $w\in A_p$. Then
there exists an $\epsilon>0$ such that $w\in A_{p-\epsilon}$. The
$\epsilon$ is not uniform and depends on $w(x)$ through its $A_p$
constant and other constants.

\bigskip
\noindent {\bf Definition:} We say $w$ satisfies a Reverse
H\"older condition of order $r>1$ if and only if for all balls $B$
one has
       $$\big({{1}\over{|B|}}\int_Bw^r(x)\, dx\big)^{1/r}\leq C
       {{1}\over{|B|}}\int_B w(x)\, dx,$$
with $C$ independent of $B$.
\medskip\noindent
{\bf Remark 3:} A pedantic remark is that the definition above is equivalent to the definition for Reverse H\"older formed by taking cubes instead
of balls $B$.

\bigskip
\noindent {\bf Theorem 3([3], [11]):} Assume $w\in A_p$. Then $w$
satisfies a Reverse H\"older condition for some $r>1$. Conversely if
$w$ satisfies a Reverse H\"older condition for some $r>1$, then $w\in
A_p$ for some $p$.

\medskip
\noindent

\medskip
\noindent {\bf Remark 4:} One also refers to $A_\infty=\cup_{p\geq
1}A_p$ and the Theorem above thus identifies $A_\infty$ with weights
that satisfy a Reverse H\"older condition.
\bigskip\noindent
{\bf Proposition 4:} Let $w\in A_p$. Then
         $$\int_{R^n} {{w(x)}\over{(1+|x|)^{np}}}\, dx\leq C<\infty.$$

\medskip\noindent
{\bf Proof:} We first claim that for any $q>p$, we have
          $$\int_{R^n} {{w(x)}\over{(1+|x|)^{nq}}}\, dx\leq C<\infty.\eqno (6)$$
By Theorem 2, since $w\in A_p$, we have $w\in A_{p-\epsilon}$ for some $\epsilon>0$. Using our claim (6) we obtain that $(6)$ holds for all $q>p-\epsilon$. Thus $(6)$ holds for $q=p$. To establish $(6)$ we use Proposition 1. We have via cutting in dyadic annuli,
 $$ \int_{R^n} {{w(x)}\over{(1+|x|)^{nq}}}\, dx\leq \int_{B_2(0)} w(x)\, dx+C\sum_{k\geq 2}2^{-knq} \int_{B_{2^k}(0)\setminus B_{2^{k-1}}(0)}w(x)\, dx,$$
where $B_{2^k}(0)$ denotes the ball of radius $2^k$ centered at the origin. The expression on the right above using Proposition 1 and the hypothesis that $w\in A_p$ is bounded by
          $$C\big(1+\sum_{k\geq 2} 2^{-kn(q-p)}\big)\int_{B_2(0)} w(x)\, dx.$$
This yields our proposition.

\medskip
We now apply the results proved above to polynomials.

\bigskip\noindent
{\bf Proposition 5:} Let $P(\xi),\ \xi\in R^n$ denote any polynomial with possibly complex valued coefficients. Then for any ball $B\subset R^n$, we have for any $r$, $1\leq r<\infty$,
             $$\bigg({{1}\over{|B|}}\int_B |P(\xi)|^r\, d\xi\bigg)^{{{1}\over{r}}}\leq C{{1}\over{|B|}}\int_B |P(\xi)|\, d\xi.\eqno(7)$$
The constant $C$ is an absolute constant which is independent of the radius and center of the ball $B$, but depends on $r,n,d$, where $d$ is the degree of $P(\xi)$.

\medskip\noindent
{\bf Remark 5:} The meaning of the proposition is that for any polynomial $P(\xi)$, $|P(\xi)|$ satisfies a Reverse H\"older condition and thus $|P(\xi)|$ is a weight in Muckenhoupt's class $A_p$ for some $p$, thanks to Theorem 3.

\bigskip\noindent
{\bf Proof:} Assume the ball $B$ has its center at $\sigma_0$ and has radius $\delta$. By translation and dilation one sees that $(7)$ is equivalent to checking,
           $$\bigg(\int_{B_1(0)} |P(\delta\xi+\sigma_0)|^r\, d\xi\bigg)^{{{1}\over{r}}}\leq C\int_{B_1(0)} |P(\delta\xi+\sigma_0)|\, d\xi,\eqno(8)$$
           where $B_1(0)$ denotes the ball of unit radius centered at the origin. Now $Q(\xi)=P(\delta\xi+\sigma_0)$ is  also a polynomial whose degree does not exceed that of $P(\xi)$. Now we use the well-known fact that any two norms on a finite dimensional vector space are equivalent. We apply this fact to the finite dimensional space of polynomials whose degree does not exceed the degree of $P(\xi)$. In particular we apply the Functional Analysis fact stated in the previous sentence to the polynomial $Q(\xi)$. We immediately obtain $(8)$, and this concludes the proof of the proposition.

\bigskip
Combining Propositions 3, 4, 5 we obtain,
\bigskip\noindent
{\bf Proposition 6:} Given any polynomial $P(\xi)$, there exists $1<p<\infty$ such that,
 $$\int_{R^n} {{|P(\xi)|^{-{{1}\over{p-1}}}}\over{(1+|\xi|)^{np^\prime}}}\, d\xi\leq C<\infty,\ {{1}\over{p}}+{{1}\over{p^\prime}}=1.$$

 \medskip
 Using the previous Proposition we have

 \bigskip\noindent
 {\bf Proposition 7:} Let $Q_0$ denote the cube of unit edge length centered at the origin. Then there exists $\xi_0\in Q_0$, such that for all lattice points $m\in Z^n$ we have
                 $$|P(\xi_0+m)|^{-{{1}\over{p-1}}}\leq C(1+|m|)^{np^\prime}.\eqno (9)$$
 The constant $C$ is independent of $m$ but may depend on $\xi_0$. The proof will show that several choices of $\xi_0$ are possible and the choices constitute a subset of full measure in $Q_0$.

 \medskip\noindent
 {\bf Proof:} We subdivide the integral in Proposition 6 over disjoint cubes using the lattice points. That is we write the conclusion of Proposition 6 in the form,
           $$ \int_{Q_0}\sum_{m\in Z^n} {{|P(\xi+m)|^{-{{1}\over{p-1}}}}\over{(1+|\xi+m|)^{np^\prime}}}\, d\xi\leq C. $$
 Thus we can find $\xi_0\in Q_0$ where the integrand is finite. This yields,
 $$ \sum_{m\in Z^n} {{|P(\xi_0+m)|^{-{{1}\over{p-1}}}}\over{(1+|\xi_0+m|)^{np^\prime}}}\leq C_1.$$
From the inequality above, our proposition is immediate by considering each individual term of the sum above.

\bigskip\noindent
{\bf Proof of Theorem 1:} We are now ready to prove our theorem. Consider solving the conjugated problem
             $$e^{-i(x,\xi_0)}P(D)(e^{i(x,\xi_0)}u)=f,\eqno (10)$$
where $\xi_0\in Q_0$. Taking the Fourier transform in $(10)$ we are left to solve the multiplier problem on the $n$-torus

             $$\hat{u}(m)={{\hat{f}(m)}\over{P(\xi_0+m)}}.$$
Using $(9)$ above we have
      $$|\hat{u}(m)|\leq C (1+|m|)^{np^\prime(p-1)}|\hat{f}(m)|.$$
Since $pp^\prime-p^\prime=p$, it follows from above that
              $$|\hat{u}(m)|^2(1+|m|^2)^\rho\leq C (1+|m|^2)^{np+\rho}|\hat{f}(m)|^2.$$
Summing in $m$, we obtain our theorem. Thus we clearly see that the estimate in our theorem is linked to the $A_p$ class in which the symbol of $P(D)$, which is a polynomial lies in. This ends the proof of the Theorem.

\bigskip\noindent
{\bf Remark 6:} We remark that the estimate in Theorem 1 is open. In the sense, using Theorem 2, we know that $|P(\xi)|$ lies in $A_{p-\epsilon}$ and thus $s$ can be improved to $n(p-\epsilon)$. Lastly we do not know if the estimate in Theorem 1 is best possible.

\bigskip\noindent
{\bf Example:} Consider the differential operator
              $$P(D)=\Pi_{j=1}^n D_j^2,\ D_j={{1}\over{i}}{{\partial}\over{\partial x_j}}.$$
The symbol is $P(\xi)=\Pi_{j=1}^n \xi_j^2$. We easily verify that $|P(\xi)|$ lies in $A_p$ with $p=3+\epsilon$ for any $\epsilon>0$.
Thus $s=(3+\epsilon)n,$ for any $\epsilon>0$. The example can easily be generalized to
           $$P(D)=\Pi_{j=1}^n D_j^{m_j},$$
with $m_j$ a natural number.

Now $s=np$ with $p>m_0+1$, where $m_0=\max\{m_1,m_2,\cdots,m_n\}$.

\vfill\eject

{\centerline{\bf REFERENCES}}

\bigskip\noindent
[1] R. Beals and C. Fefferman, On local solvability of linear partial differential equations, Annals of Math., {\bf 97}, (1973), 482-498.

\medskip\noindent
[2] S. Chanillo and F. Treves, Local exactness in a class of differential complexes, J. Amer. Math. Soc., {\bf 10}(2), (1997), 393-426.

\medskip\noindent
[3] R. Coifman and C. Fefferman, Weighted norm inequalities for maximal functions and singular integrals, Studia Math., {\bf 51}, (1974), 241-250.

\medskip\noindent
[4] P. Cordaro and J. Hounie, Local solvability for a class of differential complexes, Acta Math., {\bf 187}(2), (2001), 191-212.

\medskip\noindent
[5] L. Ehrenpreis, Solutions of some problems of division I, Amer. J. of Math., {\bf 76}, (1954), 883-903.

\medskip\noindent
[6] L. Ehrenpreis, Solutions to some problems of division II, Amer. J. of Math., {\bf 78},(1956), 685-715.
\medskip\noindent
[7] L. Ehrenpreis, Solutions to some problems of division III, Amer. J. of Math., {\bf 82}, (1960), 522-588.
\medskip\noindent
[8] C. Fefferman and B. Muckenhoupt, Two non-equivalent conditions for weight functions, Proc. Amer. Math. Soc., {\bf 45}(1), (1974), 99-104.
\medskip\noindent
[9] B. Malgrange, Existence et approximation des solutions des \'equations aux d\'eriv\'ees partielle et des \'equations de convolution, Ann. Inst. Fourier Grenoble, {\bf 6}, (1955-1956), 271-355.

\medskip\noindent
[10] B. Muckenhoupt, Weighted norm inequalities for the Hardy maximal function, Trans. Amer. Math. Soc., {\bf 165}, (1972), 207-226.

\medskip\noindent
[11] B. Muckenhoupt, The equivalence of two conditions for weight functions, Studia Math., {\bf 49}, (1974), 101-106.
\medskip\noindent
[12] L. Nirenberg and F. Treves, On local solvability of linear partial differential equations I: Necessary conditions, Comm. Pure Appl. Math., {\bf 23},(1970), 1-38.

\bigskip\noindent
Deptt. of Math.

\noindent
Rutgers University,

\noindent
110 Frelinghuysen Rd.,

\noindent
Piscataway, NJ 08854, USA

\noindent
{\tt chanillo@math.rutgers.edu}

\end